\documentclass[UTF8, a4paper, 10pt]{ctexart}
\usepackage[left=2.50cm,  right=2.50cm,  top=2.50cm,  bottom=2.50cm]{geometry} % 页边距
\CTEXsetup[format={\Large\bfseries}]{section} % 设置章标题居左

\usepackage{helvet}

\usepackage[colorlinks,
            linkcolor=blue,
            anchorcolor=blue,
            citecolor=black,
            ]{hyperref}
\usepackage{amsmath,  amsfonts,  amssymb} % math equations,  symbols
\usepackage[english]{babel}
\usepackage{color}      % color content
\usepackage{graphicx}   % import figures
\usepackage{url}        % hyperlinks
\usepackage{bm}         % bold type for equations
\usepackage{multirow}
\usepackage{booktabs}

\usepackage{epstopdf}
\usepackage{epsfig}
\usepackage{algorithm}
\usepackage{algorithmic}
\usepackage{fancyhdr} % 设置页眉,页脚
\usepackage{hyperref} %bookmarks
\hypersetup{colorlinks,  bookmarks,  unicode} %unicode

\usepackage[marginal]{footmisc}

\usepackage[marginal]{footmisc}

\title {\textbf{A note on 1-2-3 and 1-2 Conjectures for 3-regular graphs}}

\author{\normalsize Jing-zhi \textsc{Chang}$^{1}$,\quad  Chao \textsc{Yang}$^{1,*}$,\quad  Zhi-xiang \textsc{Yin}$^{1}$,\quad  Bing \textsc{Yao}$^{2}$ \\
\normalsize1. School of Mathematics, Physics and Statistics; Center of Intelligent Computing and Applied Statistics,\\
\normalsize Shanghai University of Engineering Science, Shanghai, 201620, China\\
\normalsize2. College of Mathematics and Statistics, Northwest Normal University, Lanzhou, 730070, China}

\date{} % 不显示默认日期
\begin{document}
    \maketitle
\footnote{\noindent
Supported by National Natural Science Foundation of China under Grant Nos. 61672001, 61662066, 62072296.\\
$^*$Corresponding author. E-mail address: yangchaomath0524@163.com}

% \vskip 1mm
\renewcommand{\abstractname}{} % 不显示摘要名字

\noindent   % 摘要无缩进

\noindent {\textbf{Abstract:}}
The 1-2-3 Conjecture, posed by Karo\'{n}ski, {\L}uczak and Thomason, asked whether every connected graph $G$ different from $K_2$ can be 3-edge-weighted so that every two adjacent vertices of $G$ get distinct sums of incident weights. The 1-2 Conjecture states that if vertices also receive colors and the vertex color is added to the sum of its incident edges, then adjacent vertices can be distinguished using only $\{ 1,2\}$. In this paper we confirm 1-2 Conjecture for 3-regular graphs. Meanwhile,
we show that every 3-regular graph can achieve a neighbor sum distinguishing edge coloring by using 4 colors, which answers 1-2-3 Conjecture positively.



\section{Introduction \ \ }

All considered graphs are finite, undirected, simple and
connected. Let $N(u)$ denote the set of neighbors of vertex $u$
and $0\leq s<t$. Let $d_G(v)$ be the degree of vertex $v$ of $G$.
For general theoretic notations, we follow [3].

Kar\'{o}nski et al.[7] introduced and investigated a coloring of
the edges of a graph with positive integers so that adjacent
vertices have different sums of incident edge colors. More
precisely, let $f: E \to [k]=\{1,2,\dots,k\}$ be an edge coloring
of G (such a coloring is also called a $k$-coloring of $G$). For
$x \in V$, set
 \begin{large} $$\sigma(x)=\sum_{e\ni x} f(e).$$ \end{large}
A k-coloring $f$ of $G$ is called neighbor sum distinguishing edge
coloring if $\sigma(x)\ne \sigma(y)$ for every edge $xy\in E(G)$.
In other words, the vertex coloring $\sigma$ induced by $f$ in the
above described way must be proper. The minimum integer $k$ for
which there is a neighbor sum distinguishing edge coloring of a
graph $G$ will be denoted by gndi$_{\Sigma}(G)$.

In [7] Kar\'{o}nski, {\L}uczak and Thomason posed the following
elegant problem, known as the 1-2-3 Conjecture.

\vskip 2mm

 \noindent{{\bf Conjecture 1.}}\ \ {\sl Let G be a connected graph and $G\ne
K_{2}$. Then} gndi$_{\Sigma}(G)\le 3$.

\vskip 2mm

Toward the Conjecture 1, Kalkowski, Karo\'{n}ski and Pfender [6]
showed that if $G$ is a $k$-colorable graph with $k$ odd then $G$
admits a vertex-coloring $k$-edge-weighting. So, for the class of
3-colorable graphs, including bipartite graphs, the answer is
affirmative. However, in general, this question is still open.
Addario-Berry et al.[1] showed that every graph without isolated
edges has a proper $k$-weighting when $k=30$. After improvements
to $k=15$ in [2] and $k=13$ in [10], Kalkowski, Kar\'{o}nski, and
Pfender [6] showed that every graph without isolated edges has a
proper $5$-weighting. Przybylo [9] showed that every $d$-regular
graph with $d\geq 2$ admits a vertex-coloring edge 4-weighting and
every $d$-regular graph with $d\ge 10^{8}$ admits a
vertex-coloring edge 3-weighting.

In [8] the following problem related to 1-2-3 Conjecture was
introduced. Let $f: E \cup V\to [k]=\{1,2,\dots,k\}$ be a total
k-coloring of a graph $G=(V,E)$. For every vertex $v$, let
\begin{large}$$t(v):=f(v)+\sum_{u\in N(v)}\  f(uv)=f(v)+\sigma (v),$$\end{large}
Thus, $t(v)$ is the sum of incident colors of $v$ and the color of
$v$. We say that $f$ is a neighbor total sum distinguishing total
coloring of $G$ if $t(u) \ne t(v)$ for all adjacent vertices $u$,
$v$ in $G$. Similarly as above, the minimum value of $k$ for which
there exists a total neighbor sum distinguishing coloring of a
graph $G$ will be denoted by tgndi$_{\Sigma}(G)$.

In [8] Przybylo and W\'{o}zniak posed the following problem, known
as the 1-2 Conjecture.

\vskip 2mm

\noindent {{\bf Conjecture 2.}}\ \ {\sl Let G be a connected
graph. Then} tgndi$_{\Sigma}(G)\le 2$.

\vskip 2mm

Up to now, it is known that for every graph $G$,
tgndi$_{\Sigma}(G)\le 3$ (see [5]).

This paper is organized as follows. In Section 2, we offer an
important structural lemma that every connected graph $G$ contains
a $m$-partite spanning subgraph $H$ such that
$(1-\frac{1}{m})d_G(v)\leq d_H(v)$. Naturally, every 3-regular
graph $G$ has a maximal bipartite spanning subgraph $H$ such that
$G-E(H)$ is either isolated vertices or isolated edges. Via the
structural property between $H$ and $G-E(H)$, we design a coloring
algorithm to compute gndi$_{\Sigma}(G)$ and tgndi$_{\Sigma}(G)$
for 3-regular graphs $G$ in Section 3 and Section 4, respectively. We
show that gndi$_{\Sigma}(G)\leq 4$ and tgndi$_{\Sigma}(G)\leq 2$
for any 3-regular graph $G$. Therfore, the first result answers
1-2-3 Conjecture positively and the second one confirms 1-2
Conjecture.

\section{Structural Lemma \ \ }

\noindent {{\bf  Lemma 3.}}\ \ {\sl
Let $G$ be a graph on $n$ vertices. Then it exists a $m$-partite spanning subgraph $H$ such that
$(1-\frac{1}{m})d_G(v)\leq d_H(v)$ for all $v\in V(G)$, where $m$ is a positive
integer and $m\leq n$.}

\textbf{\emph{Proof}} Let $H$ be a maximal $m$-partite spanning subgraph of $G$ with the greatest
possible number of edges. Let $\{V_1,V_2,\dots,V_m\}$ be the $m$-partition of $V(H)$ and let
$v\in V_1$, $d_{V_i}(v)=|N_{V_i}(v)|$, $N_{V_i}(v)=\{u:u\in V_i,uv\in E(G)\}$,
$i=1,2,\dots,m$. Then $d_{V_1}(v)\leq d_{V_i}(v)$, $i=1,2,\dots,m$. Otherwise, it exists an $i_0$ such that
$d_{V_1}(v)>d_{V_{i_0}}(v)$, and we use $V_1\setminus \{v\}$, $V_{i_0}\cup \{v\}$ instead of $V_1$,
$V_{i_0}$, respectively, and then it generates a new maximal $m$-partite spanning subgraph $H^{'}$ of
$G$. Obviously, $\varepsilon (H^{'})> \varepsilon (H)$, a contradiction. Therefore,

$$(m-1)d_{V_1}(v)\leq \sum_{i=2}^{m}d_{V_i}(v)=d_H(v), $$
~where~$d_H(v)=|N_H(v)|$.

Let $d_G(v)=|N_G(v)|$. Then

$$d_G(v)=d_{V_1}(v)+d_H(v)\leq \frac{1}{m-1}d_H(v)+d_H(v)=\frac{m}{m-1}d_H(v)$$

Hence
$$(1-\frac{1}{m})d_G(v)\leq d_H(v)$$ for all $v\in V(G)$. $\blacksquare$

By Lemma 3, every 3-regular graph $G$ contains a maximal bipartite
spanning subgraph $H$ such that $G-E(H)$ is either isolated
vertices or isolated edges. Via the structural between $H$ and
$G-E(H)$, we further study $gndi_{\Sigma}(G)$ and
$tgndi_{\Sigma}(G)$ of $3$-regular graphs $G$.

For the sake of
narrative, we fix some natation. Let
$G=(V_{X},V_{Y},E_{X},E_{Y},E_{H})$ be a 3-regular graph with
vertex partition $(V_{X},V_{Y})$ and edge partition
$(E_{X},E_{Y},E_{H})$, where $E_X, E_Y$ and $E_H$ represent the
edge sets in $V_X, Y_Y$ and $H$, respectively. The maximal
spanning bipartite graph $H$ is the graph with vertex set
$\{x_{i}:1\le i\le m\}\cup \{y_{i} : 1\le i\le n\}$ and edge set
$\{x_{i}y_{i}:1\le i\le m,1\le i\le n\}$. We use $a_{1}$ and
$b_{1}$ to denote the number of vertices with degree 2 and 3 in
$X$, respectively. Let $a_{2}$ and $b_{2}$ denote the number of
vertices with degree 2 and 3 in $Y$, respectively. We use
$e_{x_{i}}$ for $i=1,2,\dots ,\frac{a_{1}}{2}$ to represent the
edge with two endpoints $v_{x_{i}}$, $v'_{x_{i}}$ in $X$. Let
$e_{y_{j}}$ for $j=1,2,\dots ,\frac{a_{2}}{2}$ to represent the
edge with two endpoints $v_{y_{j}}$, $v'_{y_{j}}$ in $Y$. In the
bipartite graph $H$, let $v_{2}(x)$ and $v_{3}(x)$ be the vertex
with degree 2 and 3 in $X$, respectively. Similarly, $v_{2}(y)$
and $v_{3}(y)$ denote the vertex with degree 2 and 3 in $Y$,
respectively. The edge between $v_{2}(x)$ and $v_{2}(y)$ is
denoted by $\widetilde{e_{2}}$ , and the edge connects $v_{2}(x)$
(or $v_{3}(x)$) and $v_{3}(y)$ (or $v_{2}(y)$) is denoted by
$\widetilde{e_{2-3}}$.

\section{gndi$_{\Sigma}(G)$ for 3-regular graphs\ \ }

\noindent {{\bf Theorem 4.}}\ \ {\sl For any 3-regular graph $G$,}
gndi$_{\Sigma}(G)\leq 4$.

\textbf{\emph{Proof}} \textbf{Case 1}  $a_{1}=a_{2}=0$.

This case implies that $G$ is a 3-regular complete bipartite
graph. This result has been proved by Chang et al. [4].

\textbf{Case 2}  $a_{1}=b_{2}=0$ or $b_{1}=a_{2}=0$.

Suppose that $a_{1}=b_{2}=0$.  We color all edges in $E_{H}$ and $E_{Y}$ with 1 and 2, respectively.
Then $\sigma(v_{x})=3,\ \sigma(v_{y})=4$. To assure that $\sigma(v_{y_{j}})\ne \sigma(v'_{y_{j}})$, recolor an incident edge of $v_{y_{j}}$ (or $v'_{y_{j}}$) with 3. Then $\sigma(v_{x})$ belongs to $\{3,5,7,9\}$ and $\sigma(v_{y})$ belongs to $\{4,6\}$.

\textbf{Case 3}  $b_{1}=b_{2}=0$.

Form Lemma 3, $G$ contains a maximal bipartite spanning subgraph
$H$ such that $G-E(H)$ is either isolated vertices or isolated
edges. We color all edges in $E_{H}$ and $E_{X}$ with 1, and color
all edges in $E_{Y}$ with 2. For any edge $v_{x_{i}}v'_{x_{i}}$ in
$X$, select one edge $e_z$ from $E_H$ such that $v_{x_{i}}$ (or
$v'_{x_{i}}$) is an endpoints of $e_z$ and recolor edge $e_z$ with
3. Assume that $v_{x_{i}}$ and  $v_{y_{j}}$ are connected by
$e_{z}$ and all incident edges (except for $e_z,
v_{x_{i}}v'_{x_{i}}, v_{y_{j}}v'_{y_{j}}$) of
$v_{x_i},v'_{x_i},v_{y_j},v'_{y_j}$ keep the color 1 as before,
and we call these edges being dominated, see Fig.1. Suppose that
$v'_{x_{i}}$ and  $v_{y_{j}}$ are connected by $e_{z}$. Then
$\sigma(v_{x_{i}})=5,\ \sigma(v'_{x_{i}})=3$,
$\sigma(v_{y_{j}})=6,\ \sigma(v'_{y_{j}})=4$. Continue this
procedure $\frac{a_{1}}{2}$ times until the weights of all
adjacent vertices in $G$ are distinct.  Now we prove it
feasibility, namely, it verifies that there exists at least one
 edge in $E_{H}$ which can not be dominated after $\frac{a_{1}}{2}-1$ operations. Suppose that all edges in $H$ are dominated after $\frac{a_{1}}{2}-1$
 operations, and if there still exists a pair of adjacent vertices $v_{x_{k}}$ and $v'_{x_{k}}$ (or $v_{y_{k}}$ and $v'_{y_{k}}$) having the same weights,
 then the four incident edges of $v_{x_{k}}$ and $v'_{x_{k}}$ in $H$ must receive the same color 1. By our coloring rule, it is impossible.

 \vskip 2mm

\begin{figure}[h]
\centering
\includegraphics[height=3.0cm]{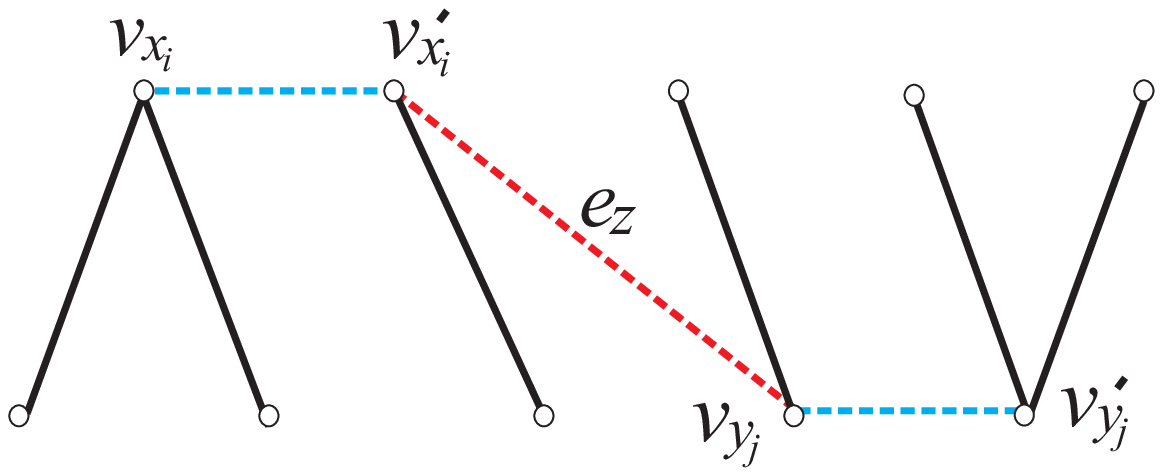}
\begin{center}
{\small Fig.1: Edges labelled by solid line are dominated.}
\end{center}
\end{figure}

\textbf{Case 4}  $a_{1}\ne a_{2},b_{1}\ne b_{2}$ and they are positive integers.

We color all edges in $E_{H}$ with 1, color all edges in $E_{X}$
with 2 and color all edges in $E_{Y}$ with 3. By the technique in
Case 3, we can distinguish all adjacent vertices whose joined by
$\widetilde{e_{2}}$. But we still need to distinguish adjacent
vertices whose joined by $\widetilde{e_{2-3}}$. Select an incident
edge from $\widetilde{e_{2-3}}$ and color it with 3, it deduces
that $\sigma(v_{x_{i}})\ne \sigma(v'_{x_{i}})$ and
$\sigma(v_{y_{j}})\ne \sigma(v'_{y_{j}})$. Moreover,
$\sigma(v_{2}(x))=4$ or 6, $\sigma(v_{2}(y))=5$ or 7,
$\sigma(v_{3}(x))$ belongs to set $\{3,5,7,9\}$ and
$\sigma(v_{3}(y))$ belongs to set $\{3,5,7,9\}$. Possibly, there
may exist some cases that the weight of adjacent vertices can not
distinguish. We deal with it as follows.

\textbf{Case 4.1} $\sigma(v_{3}(x))=3$.

\textbf{Case 4.1.1} Vertex $v_{3}(x)$ connects two vertices of degree 2 and one vertex of degree 3 in Y. Let $N(v_{3}(x))=\{v_{2}(y_{1}),v_{2}(y_{2}),v_{3}(y_{3})\}$.
Then $\sigma(v_{3}(x))=\sigma(v_{3}(y_3))=3$.

1. Vertex $v_{2}(y_{1})$ or $v_{2}(y_{2})$ has weight 7. Assume that $\sigma(v_{2}(y_{1}))=7$. Recolor $v_{3}(x)v_{2}(y_{1})$ with 2, then $\sigma(v_{3}(x))=4$, $\sigma(v_{2}(y_{1}))=8$.

2. $\sigma(v_{2}(y_{1}))=\sigma(v_{2}(y_{2}))=5$. We may recolor any incident edge of $v_{2}(y_{1})$ or $v_{2}(y_{2})$. Assume that recolor $v_{3}(x)v_{2}(y_{1})$ with 4, then $\sigma(v_{3}(x))=4$ and $\sigma(v_{2}(x_{1}))=8$.

\textbf{Case 4.1.2} Vertex $v_{3}(x)$ connects two vertices of degree 3 and a vertex of degree 2 in $Y$. Let $N(v_{3}(x))=\{v_{3}(y_{1}),v_{3}(y_{2}),v_{2}(y_{3})\}$.
Then $v_{3}(x)$ and its two neighbors all have weight 3 and $\sigma(v_{2}(y_{3}))=5$ (or 7). Recolor $v_{3}(x)v_{2}(y_{3})$ with 2 (or 4). Then $\sigma(v_{3}(x))=4$ or 6, $\sigma(v_{2}(y_{3}))=8$.

\textbf{Case 4.1.3} Vertex $v_{3}(x)$ connects three vertices of degree 3 in $Y$. Denoted by $v_{3}(y_{1})$, $v_{3}(y_{2})$ and $v_{3}(y_{3})$, respectively. If three neighbors of $v_{3}(x)$ have the same weight 3, then recolor all incident edges of $v_{3}(x)$ with 3, it follows that $\sigma(v_{3}(x))=9$ and $\sigma(v_{3}(y_{1}))=\sigma(v_{3}(y_{2}))=\sigma(v_{3}(y_{3}))=5$.

\textbf{Case 4.2} $\sigma(v_{3}(y))=5$.

 Vertex $v_{3}(x)$ connects a vertex of degree 2 with weight 7 in $Y$. Denote the vertex of degree 2 by $v_{2}(y_0)$. Recolor $v_{3}(x)v_{2}(y_0)$ with 4. Then $\sigma(v_{3}(x))=6$ and $\sigma(v_{2}(y_0))=8$.

\textbf{Case 4.3} $\sigma(v_{3}(x))=7$.

 Vertex $v_{3}(x)$ must connect two vertices of degree 2 in $Y$. Let $N(v_{3}(x))=\{v_{2}(y_{1}),v_{2}(y_{2}),v_k\}$.
 Here we need not consider the weight of $v_k$. Recolor $v_{3}(x)v_{2}(y_{1})$ and $v_{3}(x)v_{2}(y_{2})$ with 4.
 Then $\sigma(v_{3}(x))=9$ and $\sigma(v_{2}(y_{1}))=\sigma(v_{2}(y_{2}))=8$. $\blacksquare$

 \vskip 2mm

Przybylo [9] showed that every $d$-regular graph with $d\geq 2$
admits a vertex-coloring edge 4-weighting and every $d$-regular
graph with $d\ge 10^{8}$ admits a vertex-coloring edge
3-weighting. From the proof of Theorem 4, we characterize the
structural that some edges of 3-regular graphs receive color 4
when $a_1\neq a_2$ and $b_1\neq b_2$. In other words, if
Conjecture 1 is true for every 3-regular graph, then it requires a
method to deal with these edges which are colored by 4.

 \section{tgndi$_{\Sigma}(G)$ for 3-regular graphs\ \ }

\noindent {{\bf Theorem 5.}}\ \ {\sl For any 3-regular graph $G$,}
tgndi$_{\Sigma}(G)\leq 2$.

\textbf{\emph{Proof}} \ \ \textbf{Case 1}  $a_{1}=a_{2}=0$.

This case implies that $G$ is a 3-regular complete bipartite graph. We color all vertices in $X$ with $1$, color all vertices in $Y$ with 2 and color all edges in $E_{H}$ with 1. Then $t(v_{x})=4$ and $t(v_{y})=5$.

\textbf{Case 2}  $a_{1}=b_{2}=0$ or $b_{1}=a_{2}=0$.

Without loss of generality, assume that $a_{1}=b_{2}=0$.
We color all vertices in $H$ with $1$, color all edges in $E(H)$ with 1, and color all edges in $E_{Y}$ with 2.
To assure that $t(v_{y_{j}})\ne t(v'_{y_{j}})$, recolor $v_{y_{j}}$ (or $v'_{y_{j}}$) with 2. Then $t(v_{x})=4$ and $t(v_{y})$ is 5 or 6.

\textbf{Case 3}  $b_{1}=b_{2}=0$.

 We color  all vertices in $X$ with $1$, color all vertices in $Y$ with 2, color all edges in $E_{H}$ and $E_{X}$ with 1,
 and color all edges in $E_{Y}$ with 2. For any edge $v_{x_{i}}v'_{x_{i}}$ in $G$, select one edge $e_z$ from $E_H$ such that
 $v_{x_{i}}$ (or $v'_{x_{i}}$) is an endpoints of $e_z$ and recolor edge $e_z$ with 2, meanwhile, all incident edges (except for $e_z,
v_{x_{i}}v'_{x_{i}}, v_{y_{j}}v'_{y_{j}}$)
 of $v_{x_i},v'_{x_i},v_{y_j},v'_{y_j}$ keep the color 1 as before, namely they are dominated, see Fig.1. Without loss of generality, assume that $v_{x_{i}}$
 and  $v_{y_{j}}$ are connected by $e_{z}$. Then $t(v_{x_{i}})=5,\ t(v'_{x_{i}})=4$, $t(v_{y_{j}})=6,\ t(v'_{y_{j}})=7$. Continue this procedure $\frac{a_{1}}{2}$
 times until the weights of all adjacent vertices in $G$ are distinct. The feasibility of this method is similar to Case 3 in Theorem 4.

\textbf{Case 4}  $a_{1}\ne a_{2},b_{1}\ne b_{2}$ and they are positive integers.

We color all vertices in $X$ with $1$, color all vertices in $Y$ with 2, color all edges in $E_{H}$ with 2, color all edges in $E_{X}$ with 1 and color all edges in $E_{Y}$ with 2. Using the technique in Case 3 again, we change the color of $e_{z}$ from $2$ to $1$, and change the weight of vertices which are connected only by $\widetilde{e_{2}}$. By this way, we can distinguish all adjacent vertices whose are connected by $\widetilde{e_{2}}$. But we still need to distinguish adjacent vertices whose joined by $\widetilde{e_{2-3}}$. Select an incident edge from $\widetilde{e_{2-3}}$ and color it with 1, it deduces that $t(v_{x_{i}})\ne t(v^{'}_{x_{i}})$ and $t(v_{y_{i}})\ne t(v^{'}_{y_{i}})$. Meanwhile, $t(v_{2}(x))=5$ or 6, $t(v_{2}(y))=7$ or 8, $t(v_{3}(x))$ belongs to set $\{7,6,5,4\}$ and $t(v_{3}(y))$ belongs to set $\{8,7,6,5\}$. Possibly,  it still exist some cases that vertex can not distinguish. We deal it as follows.

\textbf{Case 4.1} $t(v_{3}(y))=5$.
Vertex $v_{3}(y)$ must join three vertices of degree 2 in $X$. Recolor $v_{3}(y)$ with 1, then $t(v_{3}(y))=4$ and $t(v_{2}(x))=5$.

\textbf{Case 4.2} $t(v_{3}(y))=6$.
Vertex $v_{3}(y)$ must connect two vertices of degree 2 in $X$, denoted by $v_{2}(x_{1})$ and $v_{2}(x_{2})$, respectively. Set $v_{0}\in N(v_{3}(x))-\{v_{2}(x_{1}),\ v_{2}(x_{2})\}$. Let $v'_{0}$ be the neighbor of $v_{0}$.

\textbf{Case 4.2.1} $t(v_{0})=t(v_{3}(y))=6$ and $t(v'_{0})\ne6$. Recolor $v_{3}(y)$, $v_{0}v_{3}(y)$ with 1 and $v_{0}$ with 2, then we have $t(v_{3}(y))=4$ and $t(v_{0})=6$.

\textbf{Case 4.2.2} $t(v_{0})=t(v_{3}(y))=t(v'_{0})=6$. If $d_{H}(v_{0})=2$, then $v_{0}$ connects two vertices with degree 3 in $Y$ and denoted by $v_{3}(y_{1})$ and $v_{3}(y_{2})$, respectively. Recolor $v_{0}$ with 2, it gets that $t(v_{0})=7$. If $d_{H}(v_{0})=3$, then $v_{0}$ connects two vertices of degree 3 in $Y$ and denoted by $v_{3}(y_{3})$ and $v_{3}(y_{4})$, respectively. Denote $v_{1}\in N(v_{0})-\{v_{3}(y_{3}),\ v_{3}(y_{4})\}$. Recolor $v_{0}$ and $v_{0}v_{1}$ with 2, and recolor $v_{1}$ with 1. Then $t(v_{0})=8$ and $t(v_{1})=7$.

\textbf{Case 4.3} $t(v_{3}(x))=7$. Suppose that
$N(v_{3}(x))=\{v_{y_{1}}, v_{y_{2}}, v_{y_{3}}\}$. Let
$v_{x_{1}}$, $v'_{x_{1}}$ be two neighbors of $v_{y_{1}}$ and
$v_{x_{2}}$, $v'_{x_{2}}$ be two neighbors of $v_{y_{2}}$. Then
the following cases may appear adjacent vertices having the same
weight.

\textbf{Case 4.3.1} $t(v_{y_{1}}),\ t(v_{y_{2}})$ and $t(v_{y_{3}})$ are not equal to 8. Then it exists a vertex of $N(v_{3}(x))$ having weight 7 (otherwise, $v_{3}(x)$ and its neighbors can be distinguished immediately). Recolor $v_{3}(x)$ with 2, it follows that $t(v_{3}(x))=8$.

\textbf{Case 4.3.2.} One of $\{v_{y_{1}}, v_{y_{2}}, v_{y_{3}}\}$ has weight 8 and another two vertices have weight 7.
Without loss of generality, suppose that $t(v_{y_{3}})=8$ and $t(v_{y_{1}})=t(v_{y_{2}})=7$.

\textbf{Case 4.3.2.1} $d_{H}(v_{y_{1}})=d_{H}(v_{y_{2}})=2$.

1. If $t(v_{x_{1}})=t(v_{x_{2}})=5$, recolor $v_{y_{1}}v_{3}(x)$ and $v_{y_{2}}v_{3}(x)$ with 1, then $t(v_{y_{1}})=t(v_{y_{2}})=6$, $t(v_{y_{3}})=8$ and $t(v_{x_{1}})=t(v_{x_{2}})=t(v_{3}(x))=5$.

2. Without loss of generality, assume that $t(v_{x_{1}})=6$. Recolor $v_{y_{1}}$ and $v_{y_{1}}v_{3}(x)$ with 1. Then $t(v_{y_{1}})=5$ and $t(v_{3}(x))=6$.

\textbf{Case 4.3.2.2} $d_{H}(v_{y_{1}})=2$ and $d_{H}(v_{y_{2}})=2$ hold not at the same time.
Assume that $d_{H}(v_{y_{1}})=3$ and  $t(v'_{x_{1}})=5$. Recolor $v_{y_{1}}v_{x_{1}}$, $v_{y_{1}}v_{3}(x)$ and $v_{y_{1}}$ with 1, and recolor $v_{x_{1}}$ with 2, then $t(v_{y_{1}})=4$, $t(v_{3}(x))=6$ and $t(v_{x_{1}})$ keep the same as before.

\textbf{Case 4.3.3} $\{t(v_{y_{1})},t(v_{y_{2})},t(v_{y_{3})}\}=\{6,7,8\}$. Assume that $t(v_{y_{1}})=8,\ t(v_{y_{2}})=6$ and $t(v_{y_{3}})=7$.
Recolor $v_{y_{2}}v_{3}(x)$ and $v_{y_{2}}$ with 1. Then $t(v_{y_{2}})=4$ and $t(v_{3}(x))=6$.

\textbf{Case 4.3.4} Two of $\{v_{y_{1}}, v_{y_{2}}, v_{y_{3}}\}$ have weights 8 and the remaining one has weight 7. Assume that $t(v_{y_{1}})=8,\ t(v_{y_{2}})=7$ and $t(v_{y_{3}})=8$.

1. If $d_{H}(v_{y_{2}})=2$ and $t(v_{x_{2}})=5$, recolor $v_{y_{2}}v_{3}(x)$ with 1 and recolor $v_{3}(x)$ with 2, then $t(v_{3}(x))=7$ and $t(v_{y_{2}})=6$.

2. If $d_{H}(v_{y_{2}})=2$ and $t(v_{x_{2}})=6$, recolor $v_{y_{2}}v_{3}(x)$ and $v_{y_{2}}$ with 1, then $t(v_{3}(x))=6$ and $t(v_{y_{2}})=5$.

3. If $d_{H}(v_{y_{2}})=3$, recolor $v_{y_{2}}v_{x_{2}}$, $v_{y_{2}}v_{3}(x)$ and $v_{y_{2}}$ with 1 and recolor $v_{x_{2}}$ with 2, then $t(v_{y_{2}})=4$, $t(v_{3}(x))=6$ and the weight of $v_{x_{2}}$ does not change. $\blacksquare$

\vskip 2mm

Note that Theorem 5 implies that 1-2 Conjecture is valid for
3-regular graphs. By Lemma 3, if $G$ is 4-regular, then there
exists a bipartite graph $H$ such that $2\leq d_H(u)\leq 4$ for
all $u\in V(G)$. Let $v\in V(G-E(H))$. Then $0\leq
d_{G-E(H)}(v)\leq 2$. Suppose that $X$ and $Y$ are the two
partitions of $V(H)$, then the vertex-induced subgraph of $X$ (or
$Y$) consists of cycles or isolated vertices or isolated edges or
paths. Therefore, it may offer an idea to solve 1-2-3 Conjecture
and 1-2 Conjecture for 4-regular graphs.

\section{References}
% References

\noindent {[1]}  L. Addario-Berry, K. Dalal, C. McDiarmid, B.A.
Reed, and A. Thomason, Vertex-colouring edge-weightings[J].
Combinatorica, 2007, 27: 1-12.

\noindent {[2]} L. Addario-Berry, K. Dalal, and B.A. Reed, Degree
constrained subgraphs[J]. Discrete Appl. Math. 2008, 156: 1168-1174.

\noindent {[3]} J.A. Bondy and U.S.R. Murty, Graph theory with
applications. The MaCmillan Press ltd, London and Basingstoke, New
York, 1976.

\noindent {[4]} G.J. Chang, C.H. Lu, J.J Wu, Q.L. Yu.
Vertex-coloring Edge-weightings of Graphs[J]. Taiwan. J. Math., 2011, 15(4): 1807-1813.

\noindent {[5]}  M. Kalkowski, A note on 1,2-Conjecture, to appear
in Elect. Journ. of Comb.

\noindent {[6]}  M. Kalkowski, M. Karo\'{n}ski, F. Pfender,
Vertex-coloring edge weightings: Towards the 1-2-3-conjecture[J].
J. Combin. Theory Ser. B, 2010, 100: 347-349.

\noindent {[7]} M. Karo\'{n}ski, T. {\L}uczak, A. Thomason, Edge
weights and vertex colours[J]. J. Combin. Theory ser. B, 2004, 91: 151-157

\noindent {[8]} J. Przybylo, M. Wo\'{z}niak, On a 1,2
Conjecture[J]. Discrete Math. Theor. Comput. Sci., 2010, 12(1): 101-108.

\noindent {[9]} J. Przybylo, The 1-2-3 Conjecture almost holds for
regular graphs[J]. J. Combin. Theory ser. B, 2020
https://doi.org/10.1016/j.jctb.2020.03.005.

\noindent {[10]} T. Wang and Q. Yu, On vertex-coloring
13-edge-weighting[J]. Front. Math. China, 2008, 3: 581-587.

\bibliographystyle{plain}
\renewcommand\refname{参考文献} % 参考文献标题以中文显示
\bibliography{ref}  % 参考文献

\vskip 1mm

\end{document}